\documentclass[10pt,reqno,twoside]{amsart}

\usepackage{graphicx}
\usepackage{amssymb}
\usepackage{eucal}
\usepackage{epstopdf}
\usepackage{mathrsfs}
\usepackage{amsmath}
\usepackage{amssymb}
\usepackage{enumerate}

\renewcommand{\Re}{{\mathfrak{Re}}}

\newtheorem{theorem}{Theorem}
 \newtheorem{corollary}{Corollary}[theorem]
 \newtheorem{lemma}[theorem]{Lemma}
 
 \theoremstyle{remark}

\reversemarginpar 

\begin{document}

{\def\thefootnote{}
\footnote{\today. \ {\it Mathematics Subject Classification (2000)}.
11M06, 11M26. \\
The first author is partially supported by the National Science Foundation FRG grant DMS 0244660.  The second author is supported by an NSERC research grant.
}}

\title[Lower Bounds for Moments of $\zeta'(\rho)$]
 {Lower Bounds for Moments of $\zeta'(\rho)$} 

\author{Micah B. Milinovich \and Nathan Ng}

\maketitle

\begin{abstract}
Assuming the Riemann Hypothesis, we establish lower bounds for moments of the derivative of the Riemann zeta-function averaged over the non-trivial zeros of $\zeta(s)$.  Our proof is based upon a recent method of Rudnick and Soundararajan that provides analogous bounds for moments of $L$-functions at the central point, averaged over families.
\end{abstract}

\section{Introduction}

Let $\zeta(s)$ denote the Riemann zeta-function.  In this article we are interested in obtaining lower bounds for moments of the form 
\begin{equation}\label{Jk}
J_{k}(T) =\frac{1}{N(T)} \sum_{0<\gamma\leq T} \big|\zeta'(\rho)\big|^{2k}
\end{equation}
where $k\in\mathbb{N}$ and the sum runs over the non-trivial (complex) zeros $\rho=\beta+i\gamma$ of $\zeta(s)$.  As usual, we let the function 
\begin{equation}\label{NT}
N(T) = \sum_{0<\gamma\leq T} 1 = \frac{T}{2\pi}\log\frac{T}{2\pi}-\frac{T}{2\pi} + O(\log T)
\end{equation}
denote the number of zeros of $\zeta(s)$ up to a height $T$ counted with multiplicity.  

Independently, Gonek \cite{G3} and Hejhal \cite{H} have conjectured that $J_{k}(T)\asymp (\log T)^{k(k+2)}$ for each $k\in\mathbb{R}$.  By modeling the Riemann zeta-function and its derivative using characteristic polynomials of random matrices, Hughes, Keating, and O'Connell \cite{HKO} have refined this conjecture to state that  $J_{k}(T)\sim C_{k}(\log T)^{k(k+2)}$ for a precise constant $C_{k}$ when $k\in\mathbb{C}$ and $\Re k >-3/2$.  
However, we no longer believe this conjecture to be true for $\Re k < -3/2$.  This 
is since we expect there exist infinitely many zeros $\rho$ such that $|\zeta'(\rho)|^{-1} \gg 
|\gamma|^{1/3-\varepsilon}$ for each $\varepsilon>0$.

Results of the sort suggested by these conjectures are only known for a few small values of $k$.  See, for instance, the results of Gonek \cite{G1} for the case $k=1$ and Ng \cite{N1} for the case $k=2$.  Also, Gonek \cite{G3} obtained a lower bound
in the case $k=-1$. Our main result is to obtain a lower bound for $J_{k}(T)$ for each $k\in\mathbb{N}$ of the order of magnitude that is suggested by these conjectures.

\begin{theorem}\label{th1}
Assume the Riemann Hypothesis and let $k\in \mathbb{N}$.  Then for sufficiently large $T$ we have
\begin{equation*}
\frac{1}{N(T)}\sum_{0<\gamma\leq T} \big|\zeta'(\rho)\big|^{2k} \gg_{k} (\log T)^{k(k+2)}.
\end{equation*}
\end{theorem}
Under the assumption of the Riemann Hypothesis, Milinovich \cite{M} has recently shown that $J_{k}(T)\ll_{k,\varepsilon} (\log T)^{k(k+2)+\varepsilon}$ for $k\in\mathbb{N}$ and $\varepsilon>0$ arbitrary.  When combined with Theorem \ref{th1}, this result lends strong support for the conjecture of Gonek and Hejhal for $k$ a positive integer. 

Theorem \ref{th1} can be used to exhibit large values of $\zeta'(\rho)$.  For example, as an immediate corollary we have the following result.
\begin{corollary}\label{cor}
Assume the Riemann Hypothesis and let $\rho=\tfrac{1}{2}+i\gamma$ denote a non-trivial zero of $\zeta(s)$.  Then for each $A>0$ the inequality
\begin{equation} \label{large}
  \big|\zeta'(\rho)\big| \geq (\log |\gamma|)^{A}
\end{equation}
is satisfied infinitely often.
\end{corollary}
\noindent This result was previously proven by Ng \cite{N3} by an application
of Soundararajan's resonance method \cite{S}.  The present proof is simpler and provides many more
zeros $\rho$ such that (\ref{large}) is true.  On the other hand, the resonance
method is capable of detecting much larger values of $\zeta'(\rho)$ assuming a
very weak form of the generalized Riemann hypothesis.

Our proof of Theorem \ref{th1} relies on combining a method of Rudnick and Soundararajan \cite{RS1,RS2} with a mean-value theorem of Ng (our Lemma \ref{ng}) and a well-known lemma of Gonek (our Lemma \ref{gonek}).  It is likely that our proof can be adapted to prove a lower bound for $J_{k}(T)$ of the conjectured order of magnitude for all rational $k$ (with $k\geq 1$) in a manner analogous to that suggested in \cite{RS1}.

Let $k \in\mathbb{N}$ and define, for $\xi\geq 1$, the function $\mathcal{A}_{\xi}(s) = \sum_{n\leq \xi}n^{-s}$.  Assuming the Riemann Hypothesis, we will estimate 
\begin{equation*}\label{sigma}
\Sigma_{1} = \sum_{0<\gamma\leq T}\zeta'(\rho)\mathcal{A}_{\xi}(\rho)^{k-1}\overline{\mathcal{A}_{\xi}(\rho)}^{k} \quad  \text{ and } \quad \Sigma_{2} = \sum_{0<\gamma\leq T}\big|\mathcal{A}_{\xi}(\rho)\big|^{2k}
\end{equation*}
where the sums run over the non-trivial zeros $\rho=\tfrac{1}{2}+i\gamma$ of $\zeta(s)$.  H\"{o}lder's inequality implies that
\begin{equation*}
 \sum_{0<\gamma\leq T}\big|\zeta'(\rho)\big|^{2k}\geq \frac{\big|\Sigma_{1}\big|^{2k}}{\big(\Sigma_{2}\big)^{2k-1}},
\end{equation*}
and so we see that Theorem \ref{th1} will follow from the estimates
\begin{equation}\label{bound}
 \Sigma_{1} \gg T(\log T)^{k^{2}+2} \quad \text{ and } \quad \Sigma_{2} \ll T(\log T)^{k^{2}+1}. 
\end{equation}

It is convenient to express $\Sigma_{1}$ and $\Sigma_{2}$ slightly differently. Assuming the Riemann Hypothesis, $1-\rho=\bar{\rho}$ for any non-trivial zero $\rho$ of $\zeta(s)$.  Thus, $\overline{\mathcal{A}_{\xi}(\rho)} = \mathcal{A}_{\xi}(1-\rho)$.  This allows us to re-write the sums in (\ref{sigma}) as
\begin{equation}\label{simp}
\Sigma_{1} =\!\sum_{0<\gamma\leq T}\zeta'(\rho)\mathcal{A}_{\xi}(\rho)^{k-1}\mathcal{A}_{\xi}(1\!-\!\rho)^{k} \quad  \text{ and } \quad \Sigma_{2} =\!\sum_{0<\gamma\leq T}\mathcal{A}_{\xi}(\rho)^{k}\mathcal{A}_{\xi}(1\!-\!\rho)^{k}.
\end{equation}
It is with these representations of $\Sigma_{1}$ and $\Sigma_{2}$ that we establish the bounds in (\ref{bound}).


\section{Some preliminary estimates}

For each real number $\xi\geq 1$ and each $k\in\mathbb{N}$, we define the arithmetic sequence of real numbers $\tau_{k}(n;\xi)$ by 
\begin{equation}\label{tau1}
\sum_{n\leq \xi^{k}}\frac{\tau_{k}(n;\xi)}{n^{s}} =\Big(\sum_{n\leq\xi}\frac{1}{n^{s}}\Big)^{k}= \mathcal{A}_{\xi}(s)^{k}.
\end{equation}
The function $\tau_{k}(n;\xi)$ is a truncated approximation to the arithmetic function $\tau_{k}(n)$ (the $k$-th iterated divisor function) which is defined by 
\begin{equation}\label{tau2}
\zeta^{k}(s) = \Big(\sum_{n=1}^{\infty}\frac{1}{n^{s}}\Big)^{k} = \sum_{n=1}^{\infty}\frac{\tau_{k}(n)}{n^{s}}
\end{equation}
for $\Re s> 1$.  We require a few estimates for sums involving the functions $\tau_{k}(n)$ and $\tau_{k}(n;\xi)$ in order to establish the bounds for $\Sigma_{1}$ and $\Sigma_{2}$ in (\ref{bound}).  

We use repeatedly that, for $x\geq 3$ and $k,\ell \in \mathbb{N}$, 
\begin{equation}\label{divisor1}
 \sum_{n\leq x} \frac{\tau_{k}(n)\tau_{\ell}(n)}{n} \asymp_{k,\ell} (\log x)^{k\ell}
\end{equation}
where the implied constants depend on $k$ and $\ell$. 
These bounds are well-known.  

From (\ref{tau1}) and (\ref{tau2}) we notice that $\tau_{k}(n;\xi)$ is non-negative and $\tau_{k}(n;\xi) \leq \tau_{k}(n)$ with equality holding when $n\leq \xi$.  In particular, choosing $k=\ell$ in (\ref{divisor1}) we find that, for $\xi\geq 3$, 
\begin{equation}\label{divisor2}
(\log \xi)^{k^{2}} \ll_{k} \sum_{n\leq \xi}\frac{\tau_{k}(n)^{2}}{n}\leq \sum_{n\leq \xi^{k}}\frac{\tau_{k}(n;\xi)^{2}}{n}\leq \sum_{n\leq \xi^{k}}\frac{\tau_{k}(n)^{2}}{n} \ll_{k}(\log \xi)^{k^{2}}.
\end{equation}


\section{A Lower Bound for $\Sigma_{1}$}

In order to establish a lower bound for $\Sigma_{1}$, we require a mean-value estimate for sums of the form
$$ S(X,Y;T) = \sum_{0<\gamma\leq T} \zeta'(\rho)X(\rho)Y(1-\rho) $$
where
$$ X(s) = \sum_{n\leq N} \frac{x_{n}}{n^{s}} \quad\quad \text{ and } \quad \quad Y(s) = \sum_{n\leq N}\frac{y_{n}}{n^{s}} $$
are Dirichlet polynomials.  For $X(s)$ and $Y(s)$ satisfying certain reasonable conditions, a general formula for $S(X,Y;T)$ has been established by the second author \cite{N2}.  
Before stating the formula, we first introduce some notation.  For $T$ large, we let $\mathscr{L}=\log\tfrac{T}{2\pi}$ and $N=T^{\vartheta}$ for some fixed $\vartheta\geq 0$. 
The functions $\mu(\cdot)$ and $\Lambda(\cdot)$ are used to denote the usual arithmetic functions of M\"{o}bius and von Mangoldt.  Also, we define the 
arithmetic function $\Lambda_2(\cdot)$ by $\Lambda_2(n)
=(\mu*\log^2)(n)$ for each $n\in\mathbb{N}$. 
\begin{lemma}\label{ng}
Let $x_{n}$ and $y_{n}$ satisfy $|x_{n}|,|y_{n}| \ll \tau_{\ell}(n)$ for some $\ell\in\mathbb{N}$ and assume that $0<\vartheta<1/2$.  Then for any $A>0$, any $\varepsilon>0$, and sufficiently large $T$ we have
\begin{eqnarray*}
S(X,Y;T) &=& \frac{T}{2\pi} \sum_{mn\leq N} \frac{x_{m}y_{mn}}{mn}\Big(\mathcal{P}_{2}(\mathscr{L})-2\mathcal{P}_{1}(\mathscr{L})\log n + (\Lambda*\log)(n)\Big)
\\
&& \quad - \frac{T}{4\pi} \sum_{mn\leq N} \frac{y_{m}x_{mn}}{mn}\mathcal{Q}_{2}(\mathscr{L}\!-\!\log n) + \frac{T}{2\pi} \sum_{\substack{a,b\leq N \\ (a,b)=1}}\frac{r(a;b)}{ab}\sum_{g\leq \min\big(\tfrac{N}{a},\tfrac{N}{b}\big)} \frac{y_{ag} x_{bg}}{g}
\\
&& \quad + \ O_{A}\big(T(\log T)^{-A} + T^{3/4+\vartheta/2+\varepsilon}\big)
\end{eqnarray*}
where $\mathcal{P}_1,\mathcal{P}_2$, and $\mathcal{Q}_2$ are monic polynomials of degrees 1,2,
and 2, respectively, and for $a,b\in\mathbb{N}$ the function $r(a ;b)$ satisfies the bound
\begin{equation} \label{rabbd}
  |r(a;b)| \ll  \Lambda_2(a) + (\log T) \Lambda(a) \ . 
\end{equation}
\end{lemma}
\begin{proof}
This is a special case of Theorem 1.3 of Ng \cite{N2}.
\end{proof}
Letting $\xi=T^{1/(4k)}$, we find that the choices $X(s)=\mathcal{A}_{\xi}(s)^{k-1}$ and $Y(s)=\mathcal{A}_{\xi}(s)^{k}$ satisfy the conditions of Lemma \ref{ng} with $\vartheta=1/4$, $N=\xi^{k}$, $x_{n}=\tau_{k-1}(n;\xi)$, and $y_{n}=\tau_{k}(n;\xi)$.  Consequently, for this choice of $\xi$, 
\begin{eqnarray*}
\Sigma_{1} &=& \frac{T}{2\pi} \sum_{\substack{mn\leq \xi^{k} \\ m\leq \xi^{k-1}}} \frac{\tau_{k-1}(m;\xi)\tau_{k}(mn;\xi)}{mn}\Big(\mathcal{P}_{2}(\mathscr{L})-2\mathcal{P}_{1}(\mathscr{L})\log n + (\Lambda*\log)(n)\Big)
\\
&& \quad \quad - \frac{T}{4\pi} \sum_{mn\leq \xi^{k-1}} \frac{\tau_{k}(m;\xi)\tau_{k-1}(mn;\xi)}{mn}\mathcal{Q}_{2}(\mathscr{L}\!-\!\log n)
\\
&& \quad \quad  + \frac{T}{2\pi} \sum_{\substack{a,b\leq \xi^{k} \\ (a,b)=1}}\frac{r(a;b)}{ab}\sum_{g\leq \min\big(\tfrac{N}{a},\tfrac{N}{b}\big)} \frac{\tau_{k}(ag;\xi) \tau_{k-1}(bg;\xi)}{g}+ \ O\big(T\big)
\\
&=& \mathcal{S}_{11}+\mathcal{S}_{12}+\mathcal{S}_{13} + O(T),
\end{eqnarray*}
say.  To estimate $\mathcal{S}_{11}$, notice that, for $T$ sufficiently large, $n\leq \xi^{k}=T^{1/4}$ implies that $$\Big(\mathcal{P}_{2}(\mathscr{L})-2\mathcal{P}_{1}(\mathscr{L})\log n + (\Lambda*\log)(n)\Big) \gg \mathscr{L}^{2}$$
and moreover, by (\ref{divisor2}), 
$$ \sum_{\substack{mn\leq \xi^{k} \\ m\leq \xi^{k-1}}} \frac{\tau_{k-1}(m;\xi)\tau_{k}(mn;\xi)}{mn} \geq  \sum_{n\leq \xi^{k}} \frac{\tau_{k}(n;\xi)^{2}}{n} \gg (\log T)^{k^{2}}. $$
Thus, $\mathcal{S}_{11}\gg T(\log T)^{k^{2}+2}.$ Since $\mathcal{Q}_{2}(\mathscr{L}\!-\!\log n) \ll \mathscr{L}^{2}$, we can  bound $\mathcal{S}_{12}$ by using the inequalities $\tau_{k}(n;\xi)\leq \tau_{k}(n)$ and $\tau_{k}(mn)\leq \tau_{k}(m)\tau_{k}(n)$.  In particular, by twice using (\ref{divisor1}), we find that
\begin{equation*}
\begin{split}
\mathcal{S}_{12} &\ll  T\mathscr{L}^{2} \sum_{mn\leq \xi^{k}}\frac{\tau_{k}(m)\tau_{k-1}(m)\tau_{k}(n)}{mn}
\leq  T\mathscr{L}^{2} \Bigg(\sum_{m\leq T}\frac{\tau_{k}(m)\tau_{k-1}(m)}{m}\Bigg)\Bigg(\sum_{n\leq T}\frac{\tau_{k-1}(n)}{n}\Bigg) 
\\
&\ll T(\log T)^{2+k(k-1) + k-1}\ll T(\log T)^{k^{2}+1}.
\end{split}
\end{equation*}
It remains to consider the contribution from $\mathcal{S}_{13}$. 
Again using the inequalities $\tau_{k}(n;\xi)\leq \tau_{k}(n)$ and $\tau_{k}(mn)\leq \tau_{k}(m)\tau_{k}(n)$
along with (\ref{rabbd}), it follows that $\mathcal{S}_{13}$ is bounded by 
\begin{equation*}
\begin{split}
 \sum_{a,b\leq \xi^{k}}&\frac{(\Lambda_2(a)+ (\log T)\Lambda(a))}{ab} \sum_{g\leq \xi^{k}} \frac{\tau_{k}(a)\tau_{k}(g)\tau_{k-1}(b)\tau_{k-1}(g)}{g}
\\
&\ll \sum_{a\leq T}\frac{(\Lambda_2(a)+(\log T)\Lambda(a))\tau_{k}(a)}{a}\sum_{b\leq T}\frac{\tau_{k-1}(b)}{b} \sum_{g\leq T} \frac{\tau_{k}(g)\tau_{k-1}(g)}{g}
\\
&\ll (\log T)^{2+(k-1)+k(k-1)} = (\log T)^{k^{2}+1}.
\end{split}
\end{equation*}
Combining this with our estimates for $\mathcal{S}_{11}$ and $\mathcal{S}_{12}$, we conclude that $\Sigma_{1}\gg T(\log T)^{k^{2}+2}$.


\section{An Upper Bound for $\Sigma_{2}$}

Assuming the Riemann Hypothesis, we interchange the sums in (\ref{simp}) and find that 
\begin{equation}\label{sigma2}
 \Sigma_{2} = N(T)\sum_{n\leq \xi^{k}} \frac{\tau_{k}(n;\xi)^{2}}{n} + 2\Re\sum_{m\leq \xi^{k}}\sum_{m < n\leq\xi^{k}} \frac{\tau_{k}(m;\xi)\tau_{k}(n;\xi)}{n}\sum_{0 < \gamma\leq T} \Big(\frac{n}{m}\Big)^{\rho}
 \end{equation}
where $N(T)$ denotes the number of non-trivial zeros of $\zeta(s)$ up to a height $T$.  Recalling that $\xi=T^{1/(4k)}$ and using (\ref{NT}) and (\ref{divisor2}), it follows that
\begin{equation}\label{mt}
 N(T)\sum_{n\leq \xi^{k}} \frac{\tau_{k}(n;\xi)^{2}}{n} \ll T(\log T)^{k^{2}+1}.
 \end{equation}
In order to bound the second sum on the right-hand side of (\ref{sigma2}), we require the following version of the Landau-Gonek explicit formula.

\begin{lemma}\label{gonek}
Let $x,T>1$ and let $\rho=\beta+i\gamma$ denote a non-trivial zero of $\zeta(s)$.  Then
\begin{equation*}
\begin{split}
\sum_{0 <\gamma \leq T} x^{\rho} &= -\frac{T}{2\pi}\Lambda(x) + O\big(x\log(2xT)\log\log(3x)\big)
\\
&\quad + O\Big(\log x \min\Big(T,\frac{x}{\langle x \rangle}\Big)\Big) + O\Big(\log(2T)\min\Big(T,\frac{1}{\log x}\Big)\Big) 
\end{split}
\end{equation*}
where $\langle x \rangle$ denotes the distance from $x$ to the closest prime power other than $x$ itself and $\Lambda(x)=\log p$ if $x$ is a positive integral power of a prime $p$ and $\Lambda(x)=0$ otherwise.
\end{lemma}
\begin{proof}
This is a result of Gonek \cite{G2,G4}. 
\end{proof}
\noindent Applying the lemma, we find that
\begin{equation*}
\begin{split}
\sum_{m\leq \xi^{k}}\sum_{m<n\leq\xi^{k}} \frac{\tau_{k}(m;\xi)\tau_{k}(n;\xi)}{n}\sum_{0<\gamma\leq T}\Big(\frac{n}{m}\Big)^{\rho} &=-\frac{T}{2\pi}\sum_{m\leq \xi^{k}}\sum_{m<n\leq\xi^{k}} \frac{\tau_{k}(m;\xi)\tau_{k}(n;\xi)\Lambda(\tfrac{n}{m})}{n}
\\
&\quad  +O\left(\mathscr{L}\log\mathscr{L} \sum_{m\leq \xi^{k}}\sum_{m<n\leq\xi^{k}} \frac{\tau_{k}(m;\xi)\tau_{k}(n;\xi)}{m}\right)
\\
&\quad + O\left(\sum_{m\leq \xi^{k}}\sum_{m<n\leq\xi^{k}} \frac{\tau_{k}(m;\xi)\tau_{k}(n;\xi)}{m}\frac{\log\frac{n}{m}}{\langle \frac{n}{m} \rangle}\right)
\\
&\quad  + O\left(\log T\sum_{m\leq \xi^{k}}\sum_{m<n\leq\xi^{k}} \frac{\tau_{k}(m;\xi)\tau_{k}(n;\xi)}{n \log\frac{n}{m}}\right)
\\
& = \mathcal{S}_{21}+\mathcal{S}_{22}+\mathcal{S}_{23}+\mathcal{S}_{24},
\end{split}
\end{equation*}
say.  Since we only require an upper bound for $\Sigma_{2}$ (which, by definition, is clearly positive), we can ignore the contribution from $\mathcal{S}_{21}$ because all the non-zero terms in the sum are negative.  
In what follows, we use $\varepsilon$ to denote a small positive
constant which may be different at each occurrence. 
To estimate $\mathcal{S}_{22}$, we note that $\tau_{k}(n;\xi)\leq \tau_{k}(n)\ll_{\varepsilon} n^{\varepsilon}$ which implies
$\mathcal{S}_{22} \ll T^{1/4+\varepsilon}$.  Turning to $\mathcal{S}_{23}$, we write $n$ as $qm+\ell$ with $-\frac{m}{2}<\ell\leq\frac{m}{2}$ and find that
$$ \mathcal{S}_{23} \ll  T^{\epsilon}
\sum_{m\leq\xi^{k}}\frac{1}{m} \sum_{q\leq\lfloor\frac{\xi^{k}}{m}\rfloor+1} \sum_{-\frac{m}{2}<\ell\leq\frac{m}{2}} \frac{1}{\langle q+\frac{\ell}{m} \rangle}$$  
where $\lfloor x \rfloor$ denotes the greatest integer less than or equal to $x$.  Notice that $\langle q+\frac{\ell}{m} \rangle=\frac{|\ell|}{m}$ if $q$ is a prime power and $\ell\neq0$, otherwise $\langle q+\frac{\ell}{m} \rangle$ is $\geq\frac{1}{2}$.  Hence, 
\begin{equation*}
\begin{split}
\mathcal{S}_{23} &\ll T^{\varepsilon}  \Big( \sum_{m\leq\xi^{k}}\frac{1}{m} 
  \sum_{{\begin{substack}{q\leq\lfloor\frac{\xi^{k}}{m}\rfloor+1
         \\ \Lambda(q) \ne 0}\end{substack}}}
 \sum_{1\leq\ell\leq\frac{m}{2}} \frac{m}{\ell} 
+ \sum_{m\leq\xi^{k}}\frac{1}{m} \sum_{q\leq\lfloor\frac{\xi^{k}}{m}\rfloor+1} \sum_{1\leq\ell\leq\frac{m}{2}} 1
\Big)
\\
&\ll T^{\varepsilon} 
\Big( \sum_{m\leq\xi^{k}} \sum_{q\leq\lfloor\frac{\xi^{k}}{m}\rfloor+1} 1  \Big)
\ll T^{1/4+\epsilon} \ . 
\end{split}
\end{equation*}
It remains to consider $\mathcal{S}_{24}$.  For integers $1\leq m<n\leq \xi^{k}$, let $n=m+\ell$.  Then
$$ \log\frac{n}{m} = -\log\big(1-\frac{\ell}{m}\big) > \frac{\ell}{m} . $$
Consequently, 
\begin{equation}\label{s4}
\mathcal{S}_{24} \ll T^{\epsilon}\sum_{m\leq \xi^{k}}\sum_{1 \le \ell \leq \xi^{k}} 
\frac{1}{(m+\ell) \frac{\ell}{m}} 
\ll T^{\epsilon} \xi^k =T^{1/4+\epsilon}
 . 
\end{equation}
Combining (\ref{mt}) with our estimates for $\mathcal{S}_{22}$, $\mathcal{S}_{23}$, and $\mathcal{S}_{24}$ we deduce that $\Sigma_{2} \ll
 T(\log T)^{k^{2}+1}$ which, when combined with our estimate for $\Sigma_{1}$, completes the proof of Theorem \ref{th1}.

\noindent Micah B. Milinovich\\
Math Department\\
University of Rochester\\
Rochester, NY  \\
14627 USA \\
\email{micah@math.rochester.edu}

\noindent Nathan Ng\\
Department of Mathematics and Statistics\\ 
University of Ottawa\\ 
585 King Edward Avenue\\ 
Ottawa, ON  \\
K1N 6N5 Canada \\
\email{nng@uottawa.ca}


\begin{thebibliography}{99}

\bibitem{G1} 
  S. M. Gonek, `Mean values of the Riemann zeta-function and its derivatives', \emph{Invent. Math. }75 (1984), 123-141. 

\bibitem{G2} 
  S. M. Gonek, ``A formula of Landau and mean values of $\zeta(s)$'' in \textit{Topics in Analytic Number Theory}, S. W. Graham and J. D. Vaaler, eds., (Univ. Texas Press, Austin, Tex., 1985), 92-97.

\bibitem{G3}
  S. M. Gonek, `On negative moments of the Riemann zeta-function', \emph{Mathematika }36 (1989), 71-88.

\bibitem{G4} 
  S. M. Gonek, `An explicit formula of Landau and its applications to the theory of the zeta function', \emph{Contemp. Math }143 (1993), 395-413.

\bibitem{H} 
  D. Hejhal, `On the distribution of $\log|\zeta'(1/2+it)|$', in \textit{Number Theory, Trace Formulas, and Discrete Groups}, K. E. Aubert, E. Bombieri, and D. M. Goldfeld, eds., Proceedings of the 1987 Selberg Symposium, (Academic Press, 1989), 343-370.

\bibitem{HKO} 
  C. P. Hughes, J. P. Keating, \and N. O'Connell, `Random matrix theory and the derivative of the Riemann zeta-function', \emph{Proc. Roy. Soc. London A }456 (2000), 2611-2627. 

\bibitem{M} 
  M. B. Milinovich, `Upper bounds for moments of $\zeta'(\rho)$', preprint.

\bibitem{N1} 
  N. Ng, `The fourth moment of $\zeta'(\rho)$', \emph{Duke Math J. }125 (2004) 243-266.

\bibitem{N2} 
  N. Ng, `A discrete mean value of the derivative of the Riemann zeta function', preprint.

\bibitem{N3} 
  N. Ng, `Extreme values of $\zeta'(\rho)$', preprint.

\bibitem{RS1} 
  Z. Rudnick \and K. Soundararajan, `Lower bounds for moments of L-functions', \emph{Proc. Natl. Sci. Acad. USA }102 (2005), 6837-6838.

\bibitem{RS2} 
  Z. Rudnick \and K. Soundararajan,`Lower bounds for moments of L-functions: symplectic and orthogonal examples', in \textit{Multiple Dirichlet Series, Automorphic Forms, and Analytic Number Theory},  Friedberg, Bump, Goldfeld, and Hoffstein, eds., (Proc. Symp. Pure Math., vol. 75, Amer. Math. Soc., 2006).

\bibitem{S}   K. Soundararajan, `Extreme values of L-functions', preprint.

\end{thebibliography}
\end{document}